\documentclass[11pt,letterpaper]{amsart}
\usepackage{amsmath,amsthm,amscd,amssymb,amsfonts, amsbsy}
\usepackage{latexsym, color, enumerate}
\usepackage{pxfonts}
\usepackage{mathrsfs}
\usepackage{marginnote}
\usepackage{todonotes}
\usepackage{hyperref}

\usepackage{mathtools}
\mathtoolsset{showonlyrefs}

\usepackage{thmtools}
\usepackage{thm-restate}

\usepackage{tikz}
\usetikzlibrary{snakes}
\usepackage{multicol}

\newcommand{\rvline}{\hspace*{-\arraycolsep}\vline\hspace*{-\arraycolsep}}

\usepackage{cancel}

\numberwithin{equation}{section}

\theoremstyle{plain}
\newtheorem{theorem}{Theorem}[section]
\newtheorem{lemma}[theorem]{Lemma}

\theoremstyle{definition}

\theoremstyle{remark}
\newtheorem{remark}[theorem]{Remark}

\newcommand{\bN}{\mathbb{N}}

\newcommand{\bR}{\mathbb{R}}

\newcommand\cB{\mathcal{B}}

\newcommand\cF{\mathcal{F}}

\newcommand\cL{\mathcal{L}}

\renewcommand{\vec}[1]{\boldsymbol{#1}}

\def\XXint#1#2#3{{\setbox0=\hbox{$#1{#2#3}{\int}$}
		\vcenter{\hbox{$#2#3$}}\kern-.5\wd0}}

\newcommand{\p}{\partial}
\newcommand{\epsi}{\varepsilon}
\newcommand{\dist}{\operatorname{dist}}

\begin{document}

\title[Unique continuation]{Unique continuation for Robin problems with non-smooth potentials} 

\author[Z. Li]{Zongyuan Li}
\address[Z. Li]{Department of Mathematics, City University of Hong Kong, 83 Tat Chee Avenue, Kowloon, Hong Kong SAR}
\email{zongyuan.li@cityu.edu.hk}


\begin{abstract}
	In this work, we study the unique continuation properties of Robin boundary value problems with Robin potentials $\eta \in L_{d-1+\epsi}$. Our results generalize earlier ones in which $\eta$ was assumed to be either zero (Neumann problem) or differentiable.
\end{abstract}

\maketitle

\section{Introduction}
In this paper, we study the unique continuation properties of second-order divergence form elliptic equations with Robin boundary conditions
\begin{subequations}\label{eqn-230419-0739}
	\begin{align}
		&\cL u := D_i (a_{ij}D_j u + b_i u) + W_i D_i u + Vu = 0 \,\, & \text{in}\,\,\Omega\cap B_1,\label{eqn-230419-0739-1} \\ 
		&\cB u := (a_{ij}D_j u + b_i u) \vec{n}_i = \eta u \,\, & \text{on}\,\,\p\Omega\cap B_1. \label{eqn-230419-0739-2}
	\end{align}
\end{subequations}
Here, $\vec{n}$ is the unit outward normal direction. Throughout the paper, for variable coefficient operators, we assume
\begin{equation*}
	a_{ij} = a_{ji} \quad \text{and} \quad \lambda I_d \leq (a_{ij})_{d\times d} \leq \Lambda I_d
\end{equation*}
for some $0<\lambda < \Lambda < \infty$.
Our main result reads
\begin{theorem} \label{thm-230414-0632}
	Let $\Omega \subset \bR^d$ be a $C^{1,1}$ domain, $d\geq 2$. Suppose $u \in H^1(\Omega \cap B_1)$ is non-trivial and solves \eqref{eqn-230419-0739-1}-\eqref{eqn-230419-0739-2} with coefficients satisfying $a_{ij} \in C^{0,1}(\Omega)$ and 
	\begin{equation} \label{eqn-230405-0406}
		\begin{split}
			&|b_i| + |W_i| \in L_{p, loc} (\Omega), \quad V \in L_{q, loc}(\Omega), \quad \text{and} \quad \eta \in L_{s, loc} (\p\Omega)\\
			&\text{with} \quad p > d, \quad  q > d/2, \quad \text{and} \quad s > d-1.
		\end{split}
	\end{equation}
Then we have the following.
\begin{enumerate}
	\item For every $x_0 \in \overline{\Omega}\cap B_1$, $u$ can at most vanish at a finite order, i.e., strong unique continuation property (SUCP) holds.
	\item The set $\{u=0\} \cap \p\Omega \cap B_1$ is of Hausdorff dimension not exceeding $(d-2)$. When $d=2$, for each $\rho \in (0,1)$, the set $\{u=0\} \cap \p\Omega \cap B_\rho$ contains at most finitely many points. 
\end{enumerate}
\end{theorem}
As usual, we say that $u \in H^1_{loc}$ solves \eqref{eqn-230419-0739-1}-\eqref{eqn-230419-0739-2} in the weak sense, if for every $\varphi \in C^\infty_c (\overline{\Omega} \cap B_1)$,
\begin{equation*}
	\int_{\Omega \cap B_1} - (a_{ij}D_j u + b_i u) D_i \varphi + (W_i D_i u + Vu) \varphi + \int_{\p\Omega \cap B_1} \eta u \varphi = 0.
\end{equation*}
By De Giorgi's estimate (see the end of Section \ref{sec-230725-0954}), any $H^1_{loc}$ weak solution is also H\"older continuous up to the boundary. Hence, it makes sense to discuss the Hausdorff dimension of the set $\{u=0\} \cap \p\Omega$.
\begin{remark}
	Our condition \eqref{eqn-230405-0406} is slightly scaling subcritical under $u\mapsto u_\lambda = u(\lambda\cdot)$, which we believe to be close to sharp. It would be interesting to construct the following counterexample: a non-trivial $u \in H^1_{loc}$ that locally satisfies
	\begin{equation*}
		\Delta u = 0\,\,\bR^d_+,\quad \p u/\p\vec{n} = \eta u\,\,\text{on} \,\,\p\bR^d_+
	\end{equation*} 
	with some $\eta \in L_{d-1-\epsi, loc}$. However, the set $\{u = 0\} \cap \p\bR^d_+$ has a positive surface measure.  From the discussion in \cite{MR1466583}, (SUCP) fails at every Lebesgue density point of $\{u = 0\} \cap \p\bR^d_+$.
\end{remark}
In application, the study of unique continuation properties of \eqref{eqn-230419-0739-1} - \eqref{eqn-230419-0739-2} is relevant to inverse problems that arise in non-destructive testing of materials and in modelling phenomena such as surface corrosion in metals. Typically corrosion on a surface is modeled by the Robin coefficient $\eta$ --- the ratio between the voltage ($u$) and the current ($\p u/\p\vec{n}$). Then, the goal in these applications is to recover Robin coefficients on a subset of the boundary by measuring the Cauchy data on the complementary subset. For examples of such studies, see \cite{MR3464049}, \cite{MR3681561}, and \cite{MR2745797}.
After the discussions on $\mathbb{R}^2$ with $\eta \in L_\infty$ in \cite{MR3464049}, the following question was raised in \cite[Remark~4.1]{MR3464049} for higher dimensions: 
let $u$ be harmonic in $B_1 \subset \bR^d$, $d\geq 3$, satisfying
\begin{equation} \label{eqn-230417-1152}
	0  \leq \frac{\p u}{\p \vec{n}}  \big/ u \leq C \,\,\text{on a relatively open subset}\,\, \Sigma \subset \p B_1.
\end{equation}
If the surface measure $\sigma( \{u=0\} \cap \Sigma )$ is positive,  can we conclude that $u \equiv 0$? 
By defining $\eta = \frac{\p u}{\p\vec{n}} / u$ on $\Sigma$, \eqref{eqn-230417-1152} can be rewritten in the form of \eqref{eqn-230419-0739-2} with a non-negative and bounded $\eta$. Our Theorem \ref{thm-230414-0632} (b) answers this question in a more general manner --- the Robin coefficient $\eta = \frac{\p u}{\p \vec{n}}  \big/ u$ is allowed to be sign-changing and unbounded, and the set $\{u=0\} \cap \Sigma$ is permitted to be of lower dimension.

Let us recall some history of studying unique continuation properties of second-order elliptic PDEs. Near an interior point, Koch and Tataru \cite{MR1809741} proved the (SUCP) for \eqref{eqn-230419-0739-1}, $d\geq 3$, with $b_i, W_i$ in spaces slightly smaller than $L_d$ and $V$ in spaces slightly larger than $L_{d/2}$. Similar result when $d=2$ was proved by Alessandrini in \cite{MR2834777}.  The sharpness of these requirements on the coefficients is demonstrated by several counterexamples.
At the boundary, besides (SUCP), there is also interest in Cauchy-type uniqueness problems --- specifically,  for a non-trivial solution $u$, researchers aim to determine the largest possible size of the zero Cauchy data set  $\{u = \p u/\p\vec{n} = 0\} \cap \p\Omega$. 
In the context of Dirichlet problem, i.e., \eqref{eqn-230419-0739-1} with $u=0$ on $\p \Omega \cap B_1$, this is to bound the size of $E_D := \{\p u / \p\vec{n} = 0\} \cap \p\Omega$. In \cite{MR1090434}, Lin proved (SUCP) when $b_i = 0$, $W_i, V \in L_\infty$, and $\p\Omega \in C^{1,1}$, and further demonstrated that $\dim (E_D) \leq d-2$. Lin's results were generalized to convex, $C^{1,\alpha}$, and $C^{1,Dini}$ domains in \cite{MR1363203, MR1466583, MR1415331}. More recently, Tolsa \cite{MR4544799} established for Laplace equations that if $\p\Omega$ is Lipschitz with a small constant, then $\sigma(E_D)=0$. Quantitative results were also obtained by Kenig-Zhao and Logunov et. al. in \cite{MR4444069, MR4356702}. It is conjectured that for Laplace equations, even if $\p\Omega$ is Lipschitz with a possibly large constant, $\sigma(E)=0$ should still hold.

For Neumann or Robin problems \eqref{eqn-230419-0739-1}-\eqref{eqn-230419-0739-2},  Cauchy-type uniqueness problems are concerned with the largest possible size of $E_N := \{u=0\}\cap \p\Omega$.  
Particularly for Neumann problems ($\eta = 0$), it was demonstrated in \cite{MR1466583} that (SUCP) holds on $C^{1,1}$ domains. Without much effort, this result can be extended to situations where $\eta \in C^{0,1}$, and it further proves that $\dim( E_N ) \leq d-2$. Further generalization to semilinear problems were explored in \cite{MR4104826}, allowing for $\eta \in W^{1,1}$ under a point-wise growth condition. Quantative results were also obtained for $\eta = 0$ in \cite{2021arXiv211014282B}, for $\eta \in C^{0,1}$ in \cite{2021arXiv211101766L}, and for Robin eigenfunctions (i.e., constant $\eta$) in \cite{MR4594557}.
It is worth mentioning that, in all these studies, $\eta$ is assumed to be at least differentiable. Such differentiability is crucial for deriving the major tool in these papers --- certain generalization of Almgren's monotonicity formulas.

The major novelty of the current work lies in combining two techniques to replace the use of monotonicity formulas found in earlier studies, therefore allowing us to handle rough coefficients.
The first technique involves our construction of an auxiliary function that simplifies Robin problems into homogeneous conormal problems. Our auxiliary function differs from those used in \cite{MR2745797} and \cite{MR4594557} and enables us to prove (SUCP) for non-differentiable, sign-changing $\eta$.  The second technique is a robust method for blowup analysis detailed in Section \ref{sec-230729-0745}. This method, originating from \cite{MR3952693} in the context of homogenization, has been adapted and further developed here for handling boundary value problems.
As a by product, our blowup analysis leads to the following results, which can be of independent interest.
\begin{restatable*}[]{corollary}{blowup}
	\label{cor-230725-1118}
	Let $d\geq 2$ and $\Omega \subset \bR^d$ be a $C^1$ domain with $0 \in \p\Omega$ and unit outward normal $\vec{n} = (0,\ldots,0,-1)$. Suppose that $u \in H^1(\Omega \cap B_1)$ satisfies \eqref{eqn-230419-0739-1}-\eqref{eqn-230419-0739-2} with coefficients satisfying \eqref{eqn-230405-0406} and \eqref{eqn-230727-0642}. Then, either $|u(x)| = O(|x|^N)$ for all $N>0$, or  there exists an $m \in \mathbb{N}$ and a non-trivial degree $m$ homogeneous polynomial satisfying
	\begin{equation*}
		D_i (a_{ij}(0)D_j P_m) = 0 \,\,\text{in}\,\,\bR^d_+,
		\quad
		a_{dj}(0)D_j P_m = 0 \,\, \text{on} \,\, \p \bR^d_+,
	\end{equation*}
	such that up to a subsequence,
	\begin{equation} \label{eqn-230727-0847}
		u(\lambda x) \big/ \left( \fint_{\Psi_A (\Omega \cap B_{\lambda})} |u|^2 \right)^{1/2}  \rightarrow P_m (x),
		\quad \forall x\in \bR^d_+ \cap B_1,
		\,\,\text{as} \,\, \lambda \rightarrow 0.
	\end{equation}
\end{restatable*}
Here, $\Psi_A$ is a linear transformation related to $A(0)$. See \eqref{eqn-230729-0749}. The conclusion of Corollary \ref{cor-230725-1118} seems to be new even when $\cL = \Delta$ and $\eta = 0$. The significance here is that we only require $\p\Omega \in C^1$. When $\p\Omega \in C^{1,1}$ and $\eta \in C^{0,1}$, \eqref{eqn-230727-0847} can be derived from the Almgren's monotonicity formula in \cite{MR1466583,  2021arXiv211101766L}

The rest of the paper is organized as follows. In Section \ref{sec-230725-0954}, we introduce the auxiliary function and prove (SUCP) in Theorem \ref{thm-230414-0632} (a). In Section \ref{sec-230729-0745}, we prove the key lemma in our blowup analysis --- the existence of ``limiting homogeneity'' in Lemma \ref{lem:convergence}, and provide the proof of Corollary \ref{cor-230725-1118}. With all these preparations, in Section\ref{sec-230729-0747} we prove Theorem \ref{thm-230414-0632} (b) by a dimension reduction argument.

Throughout the paper, we adopt the convention $a_i b_i = \sum_i a_i b_i$. In inequalities, the constant $C$ might vary between different lines. We use the notation $\fint_{\Omega} f : = |\Omega|^{-1} \int_{\Omega} f$ for the average, and write $W^k_p = \{u \in L_p: D^j u\in L_p, \,\,\forall j\leq k\}$ for usual Sobolev spaces.

\section{Strong unique continuation property} \label{sec-230725-0954}

In this section, we prove Theorem \ref{thm-230414-0632} (a). Since our coefficients verify conditions for the interior (SUCP) in \cite{MR2834777, MR1809741}, it suffices to prove it at boundary points.
Take
\begin{equation*}
	\widetilde{\eta}
	= 
	\begin{cases}
		\eta \,\, &\text{on} \,\, \p\Omega \cap B_1,\\
		-|\p\Omega \setminus B_1|^{-1} \int_{\p\Omega \cap B_1} \eta \,\, &\text{on} \,\, \p\Omega \setminus B_1.
	\end{cases}
\end{equation*}
Since $\widetilde{\eta} \in L_s(\p\Omega)$ with $s>1$ and $\fint_{\p\Omega} \widetilde{\eta} = 0$, by Lemma \ref{lem-230326-1020}, we can find a unique weak solution $\Psi$ to
\begin{equation*}
	D_i (a_{ij} D_j \Psi) = 0 \,\,\text{in}\,\,\Omega, \quad a_{ij}D_j \Psi \vec{n}_i = \widetilde{\eta}\,\,\text{on} \,\,\p\Omega,
\end{equation*}
and $\Psi \in W^1_{\widetilde{s}}$ for every $\widetilde{s}<ds/(d-1)$. A calculation shows that $v := u e^{{\color{blue}-}\Psi}$ solves
\begin{equation} \label{eqn-230321-0542}
	\begin{cases}
		D_i ( a_{ij} D_j v  + b_i v) + \widehat{W_i} D_i v + \widehat{V} v = 0 \quad \text{in} \,\, \Omega \cap B_1,\\
		(a_{ij} D_j v + b_i v) \vec{n}_i = 0 \quad \text{on} \,\, \p \Omega \cap B_1,
	\end{cases}
\end{equation}
where
\begin{equation*}
	\widehat{W_i} := W_i +  a_{ij} D_j \Psi, 
	\,\, 
	\widehat{V} :=  V + a_{ij} D_i\Psi D_j \Psi + b_iD_i \Psi + W_i D_i \Psi.
\end{equation*}
From \eqref{eqn-230405-0406}, $D\Psi \in L_{\widetilde{s}}$, and H\"older's inequality, we obtain $\widehat{W_i}\in L_{\widehat{p}},$ and $\widehat{V} \in L_{\widehat{q}}$, where $\widehat{p} = \min \{p ,\widetilde{s} \} > d$ and $\widehat{q} = \min\{q, \widetilde{s}/2, p\widetilde{s}/(p+\widetilde{s}) \} > d/2$ if we choose $\widetilde{s}$ to be sufficiently close to $ds/(d-1)$. Next, we apply the flattening and reflection technique in \cite{MR1466583} to the homogeneous conormal problem \eqref{eqn-230321-0542}. For this, first, from the proof of \cite[Theorem~0.8]{MR1466583} (see also \cite{2021arXiv211101766L}), there exists a $C^{1,1}$ change of variable
\begin{equation*}
	\Phi : x\in \Omega\cap B_1 \mapsto y\in \bR^d_+ \cap B_1 = \{(y', y_d): y_d > 0\}\cap B_1,
\end{equation*}
such that in new coordinates $y = \Phi (x)$, \eqref{eqn-230321-0542} becomes
\begin{equation*}
	\begin{cases}
		D_i ( \widetilde{a_{ij}} D_j \widetilde{v} \widetilde{b_{i}} \widetilde{v} ) + \widetilde{W_i} D_i \widetilde{v} + \widetilde{V} \widetilde{v} = 0 \quad \text{in} \,\, \bR^d_+ \cap B_1,\\
		(\widetilde{a_{ij}} D_j \widetilde{v} + \widetilde{b_{i}} \widetilde{v}) \vec{n}_i = 0 \quad \text{on} \,\, \p \bR^d_+ \cap B_1,
	\end{cases}
\end{equation*}
with
\begin{equation*} \label{eqn-230326-1025}
	(\widetilde{a_{ij}})_{d\times d}
	=
	\begin{pmatrix}
		\widetilde{A}_{d-1 \times d-1} & 0\\
		0 & \widetilde{a}_{dd}
	\end{pmatrix} \quad \text{on} \,\, \p\bR^d_+ \cap B_1.
\end{equation*}
That is, the conormal vector $a_{ij}\vec{n}_i$ in $x$-coordinates becomes a normal vector $\widetilde{a_{ij}} \vec{n}_i = (0,\cdots,\widetilde{a}_{dd})$ in $y$-coordinates. Moreover, we still have $\widetilde{a_{ij}} \in C^{0,1}, |\widetilde{b_i}+ \widetilde{W_i}| \in L_{\widehat{p}}, \widetilde{V} \in L_{\widehat{q}}$. 

The structure of $\widetilde{a_{ij}}$ allows us to take extensions: for $x_d <0$, let
\begin{equation*}
	\widetilde{v} (x',x_d) = \widetilde{v} (x', -x_d),
\end{equation*}
\begin{equation*}
	\widetilde{a_{ij}}(x', x_d) 
	= 
	\begin{cases}
		\widetilde{a_{ij}}(x', - x_d)\,\,\text{when either}\,\,i,j \neq d\,\,\text{or} \,\, i= j = d,\\
		- \widetilde{a_{ij}}(x', - x_d)\,\,\text{otherwise},
	\end{cases}
\end{equation*}
\begin{equation*}
	\widetilde{b}_i(x',x_d)
	=
	\begin{cases}
		\widetilde{b}_i(x', - x_d)\,\,\text{if}\,\,i\neq d,\\
		- \widetilde{b}_i(x', - x_d)\,\,\text{if}\,\,i= d,
	\end{cases}
	\quad
	\widetilde{W}_i(x',x_d)
	=
	\begin{cases}
		\widetilde{W}_i(x', - x_d)\,\,\text{if}\,\,i\neq d,\\
		- \widetilde{W}_i(x', - x_d)\,\,\text{if}\,\,i= d.
	\end{cases}
\end{equation*}
and $\widetilde{V}(x',x_d) = \widetilde{V}(x',-x_d)$.
It is not difficult to check that the extended $\widetilde{v}$ is still in $H^1$, and satisfies
\begin{equation} \label{eqn-230326-1032}
	D_i (\widetilde{a_{ij}}D_j \widetilde{v} + \widetilde{b_i} \widetilde{v}) + \widetilde{W}\cdot D \widetilde{v} + \widetilde{V} \widetilde{v} = 0, \quad \text{in}\,\, B_1
\end{equation}
with the extended coefficients still satisfying $\widetilde{a_{ij}} \in C^{0,1}$, $|\widetilde{b_i}| + |\widetilde{W_i}| \in L_{\widehat{p}}$, and $\widetilde{V} \in L_{\widehat{q}}$. Now, the strong unique continuation property follows from \cite{MR1809741} when $d\geq 3$ and \cite{MR2834777} when $d=2$. 

As a by-product of the construction above, we can also prove that any solution $u$ in Theorem \ref{thm-230414-0632} is H\"older continuous. Since $u$ can be linked with a solution $\widetilde{v}$ to \eqref{eqn-230326-1032} via a multiplication of H\"older continuous factor $e^{\Psi}$, a $C^{1,1}$ change of coordinates, and a reflection, it suffices to prove the interior H\"older continuity of $\widetilde{v}$.
	The latter follows from an interior De Giorgi estimate (cf. \cite[Theorem~8.24 and pp. 209]{MR1814364}), noting $\widetilde{p} > d$ and $\widetilde{q} > d/2$.

\section{Blowup analysis} \label{sec-230729-0745}
In this section, we discuss asymptotics for \eqref{eqn-230419-0739-1} - \eqref{eqn-230419-0739-2} at $0 \in \p\Omega$. Weaker than the requirement $a_{ij}\in C^{0,1}$ in Theorem \ref{thm-230414-0632},  in this section we only assume
\begin{equation} \label{eqn-230727-0642}
	\fint_{\Omega \cap B_r} |a_{ij} - a_{ij}(0)|^2 \rightarrow 0 \quad \text{as} \quad r\rightarrow 0.
\end{equation}
For $A(0) = (a_{ij}(0))_{d\times d}$, let $\Psi_A$ be the linear transformation with the matrix
\begin{equation} \label{eqn-230729-0749}
	\Psi_A = (\Theta A(0) \Theta^T)^{1/2} \Theta
	,\quad
	\text{where}\,\,
	\Theta = 
	\begin{pmatrix}
		I_{(d-1)\times (d-1)} & \rvline & \begin{matrix}
			-a_{1d} (0)/a_{dd} (0)\\ \cdots \\ -a_{d-1,d}(0)/a_{dd}(0)
		\end{matrix}\\
		\hline
		0 & \rvline & 1
	\end{pmatrix}
\end{equation}
Such transformation $\Psi_A$ will take $\bR^n_+$ to $\bR^n_+$, and map the conormal direction $A\vec{n}$ to normal direction.
In the following lemma, we prove that any solution to \eqref{eqn-230419-0739-1}-\eqref{eqn-230419-0739-2} with ``scaling subcritical'' coefficients either vanishs at an infinite order or has a unique limiting homogeneity.

\begin{lemma}\label{lem:convergence}
	Let $u \in H^1(\bR^d_+ \cap B_1)$ satisfy \eqref{eqn-230419-0739-1}-\eqref{eqn-230419-0739-2} on $\Omega = \bR^d_+ = \{x_d > 0\}$, with coefficients satisfying \eqref{eqn-230405-0406} and \eqref{eqn-230727-0642}.
	Then, we must have, as $r\rightarrow 0$,
	\begin{equation*}
	\begin{split}
		\log_2 \left( \fint_{\Psi_A (\bR^d_+ \cap B_{2r})} |u|^2 \big/ \fint_{\Psi_A (\bR^d_+ \cap B_{r})} |u|^2 \right)^{1/2}		&\text{either goes to infinity,}\\
		&\text{or converges to some $m \in \mathbb{N}$.}
		\end{split}
	\end{equation*}
\end{lemma}

Lemma \ref{lem:convergence} plays the role of monotonicity formulas in previous works,  and significantly relax the smoothness requirements on coefficients.  In particular, we have the (subsequence) convergence of Almgren's blowup, under very weak assumptions on $\p\Omega$ and coefficients. 
\blowup
It is worth mentioning that for every $x\in \bR^d_+ \cap B_1$, we also have $x \in \lambda^{-1}(\Omega \cap B_\lambda)$ for all small enough $\lambda$, on which $u(\lambda x)$ is defined. So, it makes sense to discuss the limit in \eqref{eqn-230727-0847}.
Meanwhile, although the limiting profile $P_m$ in \eqref{eqn-230727-0847} might depend on the particular choice of subsequences, the limiting homogeneity $m$, which comes from Lemma \ref{lem:convergence} and a flattening map, does not.

\begin{remark} 
	
	\begin{enumerate}
		\item When $a_{ij} \in C^{0,1}, b=0, |W_i| + |V| \in L_\infty$, and $\eta \in C^{0,1}$, the conclusions of Lemma \ref{lem:convergence} and Corollary \ref{cor-230725-1118} follow from the generalized Almgren's monotonicity formula in \cite{2021arXiv211101766L}. Actually, under these stronger assumptions, the limit in Lemma \ref{lem:convergence} has to be finite, i.e., (SUCP) holds.
		\item When $a_{ij} \in C^{\epsi}$ and $\Omega \in C^{1,\epsi}$, a constructive proof of Lemma \ref{lem:convergence} and Corollary \ref{cor-230725-1118} similar to \cite{MR1305956, MR1466314, 2023arXiv230710517K} is expected, which actually shows that the solution can be locally expanded as a homogeneous leading term plus higher order terms.
		However, such expansion could be wrong if we merely assume $\p\Omega \in C^1$, cf. \cite[Theorem~1.9]{2023arXiv230710517K}.
	\end{enumerate}
\end{remark}

\subsection{Proof of Lemma \ref{lem:convergence}}

The idea of Lemma \ref{lem:convergence} originated from \cite{MR3952693}, based on exploiting the rigidity (log-convexity) of the limiting harmonic functions. 
Here, we adapt the proof in \cite{MR3952693}, which works for homogenization, to our boundary value problems with scaling subcritical coefficients. 
Some other differences include: First, in \cite{MR3952693} only the boundedness of doubling indices is proved. Here, we go one step further to prove the convergence. Second, we introduce a parameter $\tau$ (see the claim below) that enables us to discuss the function limit convergence as $r\rightarrow 0$, instead of subsequence limits.

Denote $B_r^+ = \bR^d_+ \cap B_r$ and $\Gamma_r = \p \bR^d_+ \cap B_r$.
From \eqref{eqn-230729-0749}, the change of variable $y = \Psi_A x$ maps $\{x_d > 0\}$ to $\{y_d > 0\}$ with the resulted equation (in $y$ coordinates) have its leading coefficients at $0$ equal to $\delta_{ij}$ and the same integrability on lower order terms. Hence, without loss of generality, we may assume $a_{ij}(0) = \delta_{ij}$ and $\Psi_A = I_{d\times d}$.
	
	We first prove the following \textbf{Claim}: for every $\mu \notin \bN$, there exists $r_0=r_0 (u,\mu,\cL)>0$, such that for any $r<r_0$ and any $\tau \in [1/4, 1/2]$,
	\begin{equation*}
		\fint_{B_{2r}^+} |u|^2 \big/ \fint_{B_{r}^+} |u|^2 \leq 4^\mu \quad \text{implies} \quad \fint_{B_{2\tau r}^+} |u|^2 \big/ \fint_{B_{\tau r}^+} |u|^2 \leq 4^\mu.
	\end{equation*}
	We prove the claim by contradiction. Suppose not. Then there exists a sequence of $r_k \rightarrow 0$ and $\tau_k \in [1/4,1/2]$, such that
	\begin{equation} \label{eqn-230802-1032}
		\fint_{B_{2r_k}^+} |u|^2 \big/ \fint_{B_{r_k}^+} |u|^2 \leq 4^\mu ,\quad  \fint_{B_{2\tau_k r_k}^+} |u|^2 \big/ \fint_{B_{\tau_k r_k}^+} |u|^2 > 4^\mu.
	\end{equation}
	Passing to a subsequence, we may assume $\tau_k \rightarrow \tau_\infty \in [1/4, 1/2]$. Let
	\begin{equation*}
		v_k(y) := u(r_k y) \big/ (\fint_{B_{r_k}^+} |u|^2)^{1/2}.
	\end{equation*}
	Then $v_k$ satisfies
	\begin{equation*}
		\cL_k v_k = 0 \,\,\text{in} \,\,B_2^+,
		\quad
		\cB_k v_k = r_k \eta(r_k y) v_k \,\, \text{on} \,\, \Gamma_2
	\end{equation*}
	\begin{equation*}
	\text{with}\quad
		\fint_{B_1^+} |v_k|^2 = 1
		\quad\text{and}\quad
		\fint_{B_2^+} |v_k|^2 \leq 4^\mu,
	\end{equation*}
	where
	\begin{align*}
		&\cL_k v = D_i (a_{ij}(r_k y) D_j v + r_k b(r_k y) v ) + r_k W(r_k y)\cdot D v + r_k^2 V(r_k y) v,
		\\
		&\cB_k v = ( a_{ij}(r_k y)D_j v + r_k b(r_k y) ) \vec{n}_i v.
	\end{align*}
By the Caccioppoli inequality, we have $\|v_k\|_{H^1(B_{2-\epsi}^+)} \leq C\|v_k\|_{L_2 (B_2^+)} \leq C$ for every every $\epsi \in (0,1)$ and large enough $k$. Here, we have used the fact that all the exponents given in \eqref{eqn-230405-0406} are scaling subcritical.
See, for instance \cite[Proposition~2.1]{MR2944026} where an interior Caccioppoli inequality is proved when $\|b\|_{L_d} + \|W\|_{L_d} + \|V\|_{L_{d/2}}$ is small enough. Here, we can make such (scaling invariant) norms small by choosing $k$ large enough. The Robin potential $\eta$ can be treated in a similar way.
Now, by a Sobolev embedding, a trace theorem, and a diagonal argument, passing to a subsequence which is still denoted as $v_k$, we have
	\begin{equation*}
		v_k \rightarrow v_\infty \,\,
		\begin{cases}
			\text{weakly in}\,\,L_{2}(B_2^+), H^1(B_{2-\epsi}^+), L_{2^* (B_{2-\epsi}^+)},\,\,\text{and}\,\, L_{2^{**}(\Gamma_{2-\epsi})},
			\\
			\text{strongly in}\,\,L_2(B_{2-\epsi}^+) \,\,\text{and} \,\, L_2(\Gamma_{2-\epsi}),
		\end{cases}
		\forall \epsi \in (0,1).
	\end{equation*}
	Here, $2^* = 2d/(d-2)$ and $2^{**} = 2(d-1)/(d-2)$ when $d\geq 3$, and they are chosen to be sufficiently large numbers when $d=2$. Now, we can deduce from the above and \eqref{eqn-230802-1032} that
	\begin{equation} \label{eqn-230329-1215-1}
		\frac{\fint_{B_2^+} |v_\infty|^2}{\fint_{B_1^+} |v_\infty|^2} 
		\leq 
		\frac{ \liminf_k \fint_{B_2^+} |v_k|^2}{\lim_k \fint_{B_1^+} |v_k|^2}
		\leq
		4^\mu,
		\quad \text{and}\quad
		\frac{\fint_{B_{2\tau_\infty}^+} |v_\infty|^2}{\fint_{B_{\tau_\infty}^+} |v_\infty|^2} 
		=
		\frac{ \lim_k \fint_{B_{2\tau_k}^+} |v_k|^2}{ \lim_k \fint_{B_{\tau_k}^+} |v_k|^2} 
		\geq
		4^\mu.
	\end{equation}
	Now we derive the equation for $v_\infty$. Testing the $v_k$-equation by $\varphi\in C^\infty_c(B_2^+ \cup \Gamma_2)$ and using the symmetry of $a_{ij}$, after arranging terms, we have
	\begin{equation} \label{eqn-230329-1208}
		\begin{split}
			\int_{B_2^+} D_i v_k D_i \varphi
			&=
			\int_{B_2^+} \left( (\delta_{ij} -  a_{ij}(r_k y)) D_j \varphi + r_k W_i (r_k y) \varphi \right) D_i v_k
			\\&\quad +
			\int_{B_2^+} \left( - r_k b_i (r_k y) D_i \varphi + r_k^2 V(r_k y) \varphi \right) v_k
			+
			\int_{\Gamma_2} r_k \eta(r_k y) \varphi v_k.
		\end{split}
	\end{equation}
	From \eqref{eqn-230405-0406} and \eqref{eqn-230727-0642}, it is not difficult to see that
	\begin{align*}
		& \int_{B_2^+}  \left| (\delta_{ij} -  a_{ij}(r_k y))\right|^2 + \left| r_k W_i (r_k y) \right|^2 + \left| r_k b_i (r_k y)\right|^2 \,dy
		\\&\quad +
		\int_{B_2^+} \left| r_k^2 V(r_k y) \right|^{2_{*}} \,dy
		+
		\int_{\Gamma_2} | r_k \eta(r_k y) |^{2_{**}} \, d\sigma_y
		\rightarrow 0.
	\end{align*}
	Here, $2_{*}$ and $2_{**}$ are conjugate numbers of $2^*$ and $2^{**}$.
	Hence, a standard weak-strong convergence argument yields that the right-hand side of \eqref{eqn-230329-1208} converges to zero as $k$ goes to infinity. That is, $v_\infty$ satisfies
	\begin{equation*} 
		\Delta v_\infty = 0 \,\,\text{in}\,\, B_2^+, \quad \p v_\infty / \p\vec{n} = 0 \,\,\text{on} \,\, \Gamma_2.
	\end{equation*}
	Now, from \eqref{eqn-230329-1215-1} and the rigidity property in Lemma \ref{lem-230725-0930}, we must have $\fint_{B_{2r}} |v_\infty|^2 / \fint_{B_r} |v_\infty|^2 \equiv 4^k$ for some $k \in \bN$, which forces $\mu = k$. This is a contradiction since we have assumed $\mu \not\in \bN$. The claim is proved.
	
	Back to the proof of the lemma. Clearly, we only need to prove when
	\begin{equation*}
		M := \liminf_{r\rightarrow 0} \log_2 \left( \fint_{B_{2r}^+} |u|^2 \big/ \fint_{B_{r}^+} |u|^2 \right)^{1/2} < \infty.
	\end{equation*}
	In particular, $M\in [k,k+1)$ for some $k \in \bN$. Let $\epsi$ be a small number satisfying $M+\epsi < k+1$. Applying the claim iteratively with $\mu = M + \epsi$, we can easily deduce that
	\begin{equation*}
		\limsup_{r\rightarrow 0} \log_2 \left( \fint_{B_{2r}^+} |u|^2 \big/ \fint_{B_{r}^+} |u|^2 \right)^{1/2} \leq M + \epsi.
	\end{equation*}
	Since $\epsi$ can be arbitrarily small, we immediately have
	\begin{equation*}
		\lim_{r\rightarrow 0} \log_2 \left( \fint_{B_{2r}^+} |u|^2 \big/ \fint_{B_{r}^+} |u|^2 \right)^{1/2} = M.
	\end{equation*}

We are left to show that $M$ is an integer. This can be done by repeating the previous argument: consider the (subsequence) limit of $u(r\cdot)/(\fint_{B_r^+} |u|^2)^{1/2}$. The limit function is harmonic, and has the homogeneity exactly equal to $M$. This implies that $M$ must be an integer. The lemma is proved.

\subsection{Proof of Corollary \ref{cor-230725-1118}}
We first prove when $\Omega = \bR^d_+$. Now if $u$ vanishes at a finite order at $0$, the limit in Lemma \ref{lem:convergence} has  to be finite. Denote
\begin{equation*}
	u_\lambda := u(\lambda \cdot) \big/ \left( \fint_{E_\lambda} |u|^2 \right)^{1/2},
	\quad \text{where}\,\,
	E_r := \Psi_A(\bR^d_+ \cap B_\lambda).
\end{equation*}
Then $\fint_{E_1} |u_\lambda|^2 = 1$. Take a large $R$, such that $B_1^+ \subset E_R$. By Lemma \ref{lem:convergence}, we also have $\|u_\lambda\|_{L_2 (E_{2R})} \leq C$ for all small $\lambda$. Now, from the Caccioppoli inequality (again, since \eqref{eqn-230405-0406} is subcritical), we obtain $\|u_\lambda\|_{H^1 (E_{3R/2})} \leq C$. Hence, passing to a subsequence, we have $u_\lambda \rightarrow u_0$ weakly in $H^1(E_{3R/2})$ and strongly in $L_2(E_{3R/2})$. Following the argument in the proof of Lemma \ref{lem:convergence}, we can deduce that $u_0$ satisfies
\begin{equation*}
	D_i (a_{ij}(0)D_j u_0) = 0 \,\,\text{in}\,\,E_{3R/2},
	\quad
	a_{dj}(0)D_j u_0 = 0 \,\, \text{on} \,\, \p E_{3R/2} \cap \p \bR^d_+.
\end{equation*}
Hence, by separation of variable, $u_0 = \sum_j P_j(\Psi_A x)$ is a combination of homogeneous polynomials in $\Psi_A x$. Now from Lemma \ref{lem:convergence}, we must have $u_0$ itself is a degree $m$ homogeneous polynomial in $\Psi_A x$ (and hence, a homogeneous polynomial of the same degree in $x$), where $m$ is exactly the limit value given by Lemma \ref{lem:convergence}. To see the point-wise convergence, we apply a De Giorgi estimate (see the end of Section \ref{sec-230725-0954}) to obtain
\begin{equation} \label{eqn-230728-0429}
	\begin{split}
		[u_\lambda]_{C^\alpha (E_{3R/2})} 
		&=
		\lambda^\alpha  [u]_{C^\alpha (E_{3\lambda R/2})} \big/ \left( \fint_{E_\lambda} |u|^2 \right)^{1/2}
		\\&\leq
		C \lambda^\alpha (\lambda R)^{-\alpha} \left( \fint_{E_{2\lambda R}}  |u|^2 \right)^{1/2} \big/ \left( \fint_{E_\lambda} |u|^2 \right)^{1/2}
		\leq 
		C.
	\end{split}
\end{equation}
Here, we have also used (SUCP) and Lemma \ref{lem:convergence} to guarantee that for all small $\lambda$,
$$
\left( \fint_{E_{2\lambda R}}  |u|^2 \right)^{1/2} \big/ \left( \fint_{E_\lambda} |u|^2 \right)^{1/2} \leq C.
$$
Hence, pass to a further subsequence we can require the convergence $u_\lambda \rightarrow u_0$ to be also in $C^{\alpha/2} (E_R)$.

Now we deal with the general case $\Omega \in C^1$. Take a flattening map $\xi$, locally $\Omega$ to $\bR^d_+$ with $D\xi(0) = I_{d\times d}$. Now, $u\circ \xi^{-1}$ satisfies a problem locally on $\bR^d_+$, with coefficients still satisfying \eqref{eqn-230405-0406} and \eqref{eqn-230727-0642} since $D\xi \in C^0$.  From the earlier proof, there exists a subsequence $\lambda_k \rightarrow 0$, such that
	\begin{equation*}
	u\circ\xi^{-1} (\lambda_k y) \big/ \left( \fint_{\Psi_A (\bR^d_+ \cap B_{\lambda_k})} |u\circ \xi^{-1}|^2 \right)^{1/2}  \rightarrow P_m (y),
	\quad \forall y \in \bR^d_+ \cap B_1.
\end{equation*}
Finally, note that $\|u\circ\xi^{-1} - u\|_{L_\infty (B_r \cap \Omega\cap\bR^d_+)} \rightarrow 0$ as $r\rightarrow 0$ by the H\"older continuity of $u$. From $\xi \in C^1$ with $D\xi(0) = I_{d\times d}$ and a similar reasoning, we obtain
\begin{equation*}
	\fint_{\Psi_A (\bR^d_+ \cap B_{\lambda_k})} |u\circ \xi^{-1}|^2 \big/ \fint_{\Psi_A (\Omega \cap B_{\lambda_k})} |u|^2 \rightarrow 0.
\end{equation*}
Combining these, we reach \eqref{eqn-230727-0847}.

\section{Dimension reduction} \label{sec-230729-0747}

In this section, we prove Theorem \ref{thm-230414-0632} (b) by Federer's dimension reduction argument. By the flattening map in the reflection technique in \cite{MR1466583} (see Section \ref{sec-230725-0954}), without loss of generality, we may assume $\Omega = \bR^d_+$.

Let $\mathcal{F}$ be the set formed by all the subsets $E \subset \Gamma_1 (=\p \bR^d_+ \cap B_1)$, verifying the following: there exists a function $w \in H^1(\bR^d_+ \cap B_1)$ which is not constantly zero and solves \eqref{eqn-230419-0739-1}-\eqref{eqn-230419-0739-2} with $\Omega = \bR^d_+$ and coefficients satisfying $a_{ij}\in C^{0,1}$ and \eqref{eqn-230405-0406}, such that
\begin{equation*}
	E = \{x\in \Gamma_1: w(x) = 0\}.
\end{equation*}
We call such $w$ a defining function of $E$.  By a De Giorgi estimate, $w \in C^\alpha (B_1^+ \cup \Gamma_1)$. Hence, $E$ is always well-defined and is a relative closed subset of $\Gamma_1$. Meanwhile, $E$ cannot be the whole $\Gamma_1$ by the simple unique continuation property (since $w=0$ on $\Gamma_1$ implies $\p w/\p\vec{n} = 0$ on $\Gamma_1$ as well, from the boundary condition). 
We verify the following two properties of $\cF$.

\noindent\textbf{P.1. (Closure under appropriate scaling and translation)}: If $\lambda \in (0,1), y \in\Gamma_{1-\lambda},$ and $E\in \cF$, then $E_{y,\lambda}(:= \lambda^{-1}(E-y)) \in \cF$.

\noindent\textbf{P.2. (Existence of homogeneous degree zero ``tangent set'')}: For $y \in \Gamma_1$, $\{\lambda_k\} \downarrow 0$, and $E \in \cF$, there exists a subsequence $\{\lambda_{k'}\}$ and $F \in \cF$, such that $E_{y,\lambda_{k'}} \rightharpoonup F$ and $F_{0,\lambda} = F$ for each $\lambda>0$. Here, for $E^{(k)}, E^{(\infty)} \in \cF$, the convergence $E^{(k)} \rightharpoonup E^{(\infty)}$ means
\begin{equation*}
	\forall \epsi>0, \exists k(\epsi)\,\,\text{such that}\,\, \Gamma_{1-\epsi}\cap E^{(k)} \subset \{ x\in \Gamma_1: \dist(E^{(\infty)},x) <\epsi \}, \quad \forall k > k(\epsi).
\end{equation*}

Once \textbf{P.1} and \textbf{P.2} are verified, the conclusion in Theorem \ref{thm-230414-0632} (b) follows immediately from
\begin{lemma}[Theorem~A.4 in \cite{MR756417}, Lemma~2.2 in \cite{MR1090434}]
	Let $\cF$ be a collection of relatively closed proper subsets of the unit ball in $\bR^{d-1}$ that verifies \textbf{P.1} and \textbf{P.2}. Then, the Hausdorff dimension of each set in $\cF$ is no greater than $d-2$. When $d=2$, $E\cap [-\rho,\rho] (\subset \bR^1)$ is finite for every $E\in \cF$ and every $\rho<1$.
\end{lemma}

We are left to verify \textbf{P.1} and \textbf{P.2}.

Verifying \textbf{P.1.} is straightforward. For the set $E$ given above, let $w$ be a defining function. By (SUCP) in Theorem \ref{thm-230414-0632} (a), $w$ cannot be identically zero on any open subset of $B_1^+$. Hence, $w_{y,\lambda}(\cdot) := w(y+\lambda\cdot) \not\equiv 0$.  Clearly, $w_{y,\lambda} \in H^1$ satisfies 
\begin{equation*} 
	\mathcal{L}_{y, \lambda} w_{y,\lambda} = 0 \,\,\text{in}\,\,B_2^+, \quad \mathcal{B}_{y,\lambda} w_{y, \lambda} = \lambda \eta(y + \lambda\cdot) w_{y,\lambda}
\end{equation*}
where
\begin{equation*}
	\mathcal{L}_{y,\lambda} v = D_i (a_{ij}(y + \lambda\cdot) D_j v  + \lambda b(y + \lambda\cdot) v ) + \lambda W(y + \lambda\cdot)\cdot D v + \lambda^2 V(y + \lambda\cdot)  v
\end{equation*}
and $\mathcal{B}_{y,\lambda}$ is the associated conormal operator. It is not difficult to see that coefficients of $\mathcal{L}_{y,\lambda}$ still satisfy $a_{ij}(y + \lambda\cdot) \in C^{0,1}$ and \eqref{eqn-230405-0406}. Hence, $w_{y,\lambda}$ is a legitimate defining function for $E_{y,\lambda}$. So, $E_{y,\lambda} \in \cF$. This proves \textbf{P.1.}.

Verifying \textbf{P.2} requires more work, in which Lemma \ref{lem:convergence} is used.
Let $w$ be a defining function of $E$. By the reasoning in verifying \textbf{P.1}, we have $E_{y,\lambda} \in \cF$ with $w_{y,\lambda}$ being a defining function. By a linear transformation, without loss of generality, we may assume $a_{ij}(y) = \delta_{ij}$, and hence, $\Psi_A = I_{d\times d}$ in Lemma \ref{lem:convergence}.
Define $\widetilde{w_{y,\lambda}} := w_{y,\lambda} \big/ \left(\fint_{B_{2\lambda}^+(y)} |w_{y,\lambda}|^2\right)^{1/2}$, from which clearly $\fint_{B_2^+} |\widetilde{w_{y,\lambda}}|^2 = 1$. From \eqref{eqn-230728-0429}, we have $\| \widetilde{w_{y,\lambda}} \|_{C^\alpha(B_{3/2}^+ \cup \Gamma_{3/2})} \leq C$ for all small enough $\lambda$.
Hence, passing to a subsequence $\lambda_{k'}$, we have $\widetilde{w_{y,\lambda_{k'}}} \rightarrow w^{(\infty)}$ in $C^{\alpha/2}(B_1^+ \cup \Gamma_1)$. So,
\begin{equation*}
	E_{y,\lambda} \rightharpoonup F:= \{x \in \Gamma_1: w^{(\infty)} (x) = 0\}.
\end{equation*}
We are left to check $F_{0,\lambda} = F$ for each $\lambda>0$, i.e., $F$ is homogeneous. This follows from the homogeneity of $w^{(\infty)}$, which in turn is a consequence of Corollary \ref{cor-230725-1118}. Hence, \textbf{P.2} is verified.

With all above, Theorem \ref{thm-230414-0632} is proved.

\appendix
\section{Monotonicity and rigidity property of Neumann problems}
For a function $v \in H^1$, define the doubling index and the (averaged) Almgren's frequency function as
\begin{equation*}
	N (r) := \fint_{B_{2r}^+} |v|^2 \big/ \fint_{B_r^+} |v|^2
	\quad \text{and} \quad
	F (r) := \fint_{B_r^+} |\nabla v|^2 (r^2 - |x|^2) \big/ \fint_{B_r^+} |v|^2.
\end{equation*}
\begin{lemma} \label{lem-230725-0930}
		Let $d\geq 2$ and $v$ be a non-trivial solution to 
		\begin{equation} \label{eqn-230315-1128}
			\begin{cases}
				\Delta v = 0 \quad &\text{in}\,\,B_2^+,\\
				\p v /\p\vec{n} = 0 \quad &\text{on}\,\,\Gamma_2.
			\end{cases}
		\end{equation}
		Then both $F$ and $N$ are non-decreasing for $r \in (0,2)$. Moreover, if either $F (t) = F(s)$ or $N(t) = N(s)$ for some $t>s$, then $v$ must be homogeneous of degree $m (\in \bN)$, $F \equiv 2m$, and $N \equiv 4^{m}$.
\end{lemma}

Lemma \ref{lem-230725-0930} can be proved by an even extension and the interior result (cf., \cite{MR3952693}).
Here, we include a proof for completeness.

\begin{proof}
By calculation, we obtain
\begin{equation} \label{eqn-240606-0646}
\frac{d}{dr} \fint_{B_r^+} |v|^2 
=
2r^{-1} \fint_{B_r^+} v (x\cdot \nabla v) 
=
r^{-1} \fint_{B_r^+} |\nabla v|^2 (r^2 - |x|^2).
\end{equation}
In other words,
\begin{equation}\label{eqn-240614-1027}
\frac{d}{d\log r} \left( \log \fint_{B_r^+} |v|^2 \right)  = F(r).
\end{equation}
Also,
\begin{equation} \label{eqn-230727-1021}
	\frac{d}{dr}F(r) = 4 (r \fint_{B_r^+} |v|^2)^{-1} \left( \left(\fint_{B_r^+} |v|^2\right) \left( \fint_{B_r^+} |x\cdot \nabla v|^2 \right) -  \left( \fint_{B_r^+} v (x\cdot \nabla v) \right)^2 \right).
\end{equation}
The calculations for \eqref{eqn-240606-0646} and \eqref{eqn-230727-1021} are standard, based on integration by parts, \eqref{eqn-230315-1128}, and $x\cdot \vec{n} = 0$ on $\Gamma_2$. Details can be found in, for instance, \cite{MR1466583, MR4104826, MR0833393}.

From \eqref{eqn-230727-1021} and H\"older's inequality, $F' \geq 0$. Now, suppose $F(t) = F(s)$ for some $t>s$. This implies $F'(r) = 0$ for any $r \in (s,t)$. Hence, by the condition for achieving ``='' in H\"older's inequality, we must have $x\cdot \nabla v = C_r v$ for any $|x|=r$ and any $r \in (s,t)$. Taking the spherical harmonic expansion $v = \sum_k a_k r^k \psi_k(\theta)$,  this implies $\sum_k (C_r - k) a_k r^k \psi_k(\theta) = 0$ which is only possible when merely one of $a_k$ is non-zero, i.e. $v$ itself is homogeneous.
Let $m (\in \mathbb{N})$ be the homogeneity of $v$. Clearly, $N \equiv 4^m$. Moreover, combining $x\cdot \nabla v \equiv m v$, which comes from the homogeneity, and the second equality in \eqref{eqn-240606-0646}, we obtain $F \equiv 2m$.

Baed on the discussion on $F$ above, the corresponding results for $N$ follows from the following formula
\begin{align*}
N(r)
=
\exp{\int_r^{2r} \frac{d}{ds} \log (\fint_{B_s^+} |v|^2) \,ds}
=
\exp{\int_{\log_2 r}^{1+\log_2 r} F(2^s) \,ds},
\end{align*}
which can be derived from \eqref{eqn-240614-1027}.
\end{proof}

\section{A solvability result of Neumann problems}
In this section, we state a solvability result for
\begin{equation} \label{eqn-230405-0339}
	\begin{cases}
		D_i (a_{ij} D_j u ) = 0 \quad \text{in} \,\, \Omega,\\
		a_{ij} D_j u \vec{n}_i = g \quad \text{on} \,\, \p \Omega.
	\end{cases}
\end{equation}
\begin{lemma} \label{lem-230326-1020}
	Suppose $\Omega \subset \bR^d$ is a bounded $C^{1,1}$ domain, $a_{ij} \in C^{0,1}$, and $p \in (1,\infty)$. Then for every $g\in L_p (\p\Omega)$ with $\fint_{\p\Omega} g = 0$, there exists a unique weak solution $u$ to \eqref{eqn-230405-0339}. Moreover, $u \in W^1_q (\Omega)$ for every $q < dp/(d-1)$.
\end{lemma}
Lemma \ref{lem-230326-1020} is a standard elliptic regularity result, for which we only sketch the proof. Note that since $a_{ij} \in C^{0,1}$, \eqref{eqn-230405-0339} has both divergence and non-divergence structures: $D_i (a_{ij} D_j u ) = a_{ij}D_i D_j u + D_i(a_{ij}) D_j u$.  Hence, noting $\p\Omega \in C^{1,1}$, we have both $W^2_p$ and $W^1_p$ estimates for every $p \in (1,\infty)$, with $g$ coming from corresponding trace spaces. By interpolation, we reach the desired a priori estimate when $g \in L_p(\p\Omega)$. The solvability follows from a standard density argument: For a sequence of $g_k \in C^\alpha(\p\Omega)$ with $g_k \rightarrow g$ in $L_p(\p\Omega)$, solve \eqref{eqn-230405-0339} for $u_k \in C^{1,\alpha}(\Omega)$ by standard Schauder theory. The convergence of $u_k$ to a solution $u$ and the regularity of $u$ follow from the aforementioned a priori estimate.

\section*{Acknowledgment}
The author wishes to thank Hongjie Dong, Dennis Kriventsov, and Fanghua Lin for several inspiring discussions. The author was partially supported by an AMS-Simons travel fund. The paper was finished while the author was at Rutgers university. We would like to express sincere thanks for the anonymous reviewers whose comments/suggestions helped improve and clarify this manuscript.
%

%

\end{document}